\documentclass[12pt]{article}
\usepackage{amsmath,amssymb,amsthm,amsfonts,amstext,mathtools,dsfont,
	fullpage,xcolor,indentfirst,centernot}
\usepackage[english]{babel}
\usepackage{graphicx}

\newtheorem{theorem}{Theorem}[section]
\newtheorem{proposition}{Proposition}[section]
\newtheorem{example}{Example}[section]
\newtheorem{remark}{Remark}[section]

\begin{document}
	
\title{On convergence of Baum-Katz series for elements of linear autoregression \thanks{Supported by the grant 0118U003614 from
		Ministry of Education and Science of Ukraine (project N 2105$\Phi$).}}
	
	\author{Maryna Ilienko\\
		{\small\it  National Technical University of Ukraine 
		\small\it	``Igor Sikorsky Kyiv Polytechnic Institute''}\\
		{\small\it mari-run@ukr.net}}
	
	\date{}

	\maketitle

\begin{abstract}
We study complete convergence and closely related Hsu-Robbins-Erd\H{o}s-Spitzer-Baum-Katz series for sums whose terms are elements of linear autoregression sequences. We obtain criterions for convergence of this series expressed in moment assumptions, which for \lq\lq weakly dependent\rq\rq\ sequences are the same as in classical results concerning independent case.       


\end{abstract}

\section{Introduction}\label{intro}
Let $(X_{n},n\geq1)$ be a sequence of independent copies of a random variable (r.v.) $X$, and $S_n=\sum_{k=1}^n X_k,$ $n\geq1.$
The concept of complete convergence was introduced by Hsu and Robbins, \cite{Hsu}, and reads as follows. \emph{A random sequence $(U_n, \  \ n\geq 1)$ completely converges to a constant $C$, if}
\begin{equation*}
\sum_{n=1}^\infty\mathbb{P} \{|U_n-C|>\varepsilon  \} <\infty, \ \ \emph{for any} \ \  \varepsilon>0.
\end{equation*}
This, particularly, implies that $U_n\underset{n\rightarrow \infty} \longrightarrow C $ almost surely.
In their paper, Hsu and Robbins proved the sufficient part of the following proposition, while the converse was provided two years later by Erd\H{o}s.
\begin{proposition}[Hsu-Robbins-Erd\H{o}s, \cite{Hsu,Erd}]\label{HsuR}
	For any $\varepsilon>0$,
	\begin{equation*}
	\sum_{n=1}^\infty\mathbb{P} \Big\{\Big|\frac{S_n}{n}-\mathbb EX\Big|>\varepsilon  \Big\} <\infty \quad \text{iff} \quad \mathbb E X^{2}<\infty.
	\end{equation*} 
\end{proposition}
Closely related to Hsu-Robbins-Erd\H{o}s result are the following two no less celebrated results by Spitzer and by Baum and Katz. 

\begin{proposition}[Spitzer, \cite{Sp}]\label{SP}
	For any $\varepsilon>0$,
	\begin{equation*}
	\sum_{n=1}^\infty \frac{1}{n}\mathbb{P} \Big\{\Big|\frac{S_n}{n}-\mathbb EX\Big|>\varepsilon  \Big\} <\infty  \quad \text{iff} \quad \mathbb E|X|<\infty.
	\end{equation*}
\end{proposition}

\begin{proposition}[Baum-Katz, \cite{BK}]\label{BaumK}
	Let $0<p<2$ and $r\geq p$. Then for any $\varepsilon>0$,
	\begin{equation*}
	\sum_{n=1}^\infty n^{\frac{r}{p}-2}\mathbb{P} \Big\{\frac{|S_n|}{n^{1\!/\!p}}>\varepsilon  \Big\} <\infty  \quad \text{iff} \quad \mathbb E|X|^{r}<\infty,
	\end{equation*}
	where $\mathbb E X=0$, when  $r\geq1.$
	
\end{proposition}

Obviously, Proposition~{\rm{\ref{BaumK}}} on the convergence of Baum-Katz series covers both Hsu-Robbins-Erd\H{o}s and Spitzer results with $r=2p=2$ and $r=p=1$ respectively.

Propositions~{\rm{\ref{HsuR}}}--{\rm{\ref{BaumK}}} are fundamental facts of probability theory and have been extended in several directions. Among these extensions we distinguish results concerning complete convergence and convergence of Baum-Katz series for weighted sums of independent r.v.'s, also known as rowwise independent random arrays, (see, for instance, \cite{HuMoritz,gut,Hu,Cai,ChenMaSung} and references therein), as well as some  dependent patterns (see, for instance, \cite{amini,Balka} and references therein). 

In this paper, we focus on schemes with linear dependence, namely linear autoregression sequences, and aim to obtain analogues of Propositions~{\rm{\ref{HsuR}}}--{\rm{\ref{BaumK}}} for such sequences. In the case of \lq\lq weak dependence\rq\rq\ it is natural to expect  that the Baum-Katz series will converge under moment assumptions, similar to those in the independent setting.

On a common probability space $(\Omega, \mathcal{F},\mathbb{P})$ consider a linear autoregression sequence $(\xi_{k}, k\geq1)$, described by the following system of recurrence equations:
\begin{equation}\label{model}
\xi_{1}=\theta_{1}, \ \  \xi_{k}=q_k\xi_{k-1}+\theta_{k}, \ \ k\ge 2,
\end{equation}
where $(q_k)$ is a sequence of reals, and $(\theta_{k})$ is a sequence of independent r.v.'s. For more details concerning model \eqref{model} and its generalizations, as well as some applications, see, for instance, \cite{mikosh}, and numerous references therein.
Set
\begin{equation*}
S_n=\sum_{k=1}^n \xi_k, \ \ n\geq1.
\end{equation*}

For sequences of type \eqref{model} and some their extensions, we previously studied assumptions providing almost sure convergence of series $\sum_{k=1}^{\infty} \xi_k$ as well as series $\sum_{n=1}^{\infty} \dfrac{S_n}{n^{1+1\!/\!r}}$, $r>0$ (see, for instance, \cite{runovska1,buld,ilienko2016,ilienko2017}).

In this paper, we study necessary and sufficient conditions for the convergence of Baum-Katz series 
\begin{equation}\label{BKseries}
\sum_{n=1}^\infty n^{\frac{r}{p}-2}\mathbb{P} \Big\{\frac{|S_n|}{n^{1\!/\!p}}>\varepsilon  \Big\},
\end{equation}
where  $0<p<2$ and $r\geq p$.

Set
\begin{equation}\label{weights}
a(n,k)=\begin{cases}
0, & 1\leq n<k; \\
1, & n=k; \\
1+\sum_{l=1}^{n-k}\Big(\prod_{j=k+1}^{k+l}q_{j}\Big), & n>k. \end{cases}
\end{equation}
It is easily seen that elements of the sequence of partial sums $(S_n)$ may be represented in the following form: 
\begin{equation}\label{partial sums}
S_n=\sum_{k=1}^n a(n,k)\theta_k, \ \ n \geq 1.
\end{equation}

Representation \rm{\eqref{partial sums}} means that \textit{sums of elements of autoregression sequences can be treated as weighted sums of independent r. v.'s} $(\theta_k)$. This approach permits us to follow some ideas developed so far for weighted sums of independent r. v.'s. To be specific, we rely upon \cite{gut}, borrowing some tools to obtain necessary and sufficient assumptions providing convergence of series (\ref{BKseries}) for any $r\geq p$. 

Note that in \cite{gut} for rowwise independent random arrays sufficient conditions providing convergence of Baum-Katz series are considered only for $r=2p$ and $r=p$. Therefore, in these two specific cases due to representation \rm{\eqref{partial sums}}, similar results for autoregression sequences are immediate from \cite{gut}.
Emphasize also, that our goal is to obtain sufficient as well as necessary assumptions for series (\ref{BKseries}) to converge for any $r\geq p$.

\section{Main results}\label{Main results}

In what follows, we consider linear autoregression model \rm{\eqref{model}}, where  $q_k=q=const$,  $|q|\leq1$, and $(\theta_{k})$ is a sequence of independent copies of a random variable $\theta$. Observe that for  $-1\leq q<1$, by (\ref{weights})
\begin{equation*}	
a(n,k)=\begin{cases}
0, & 1\leq n<k; \\
1, & n=k; \\
\dfrac{1-q^{n-k+1}}{1-q}, & n>k, \end{cases}
\end{equation*}
and the triangular array $\big(a(n,k), 1\leq k\leq n, n\geq 1 \big)$ is uniformly bounded, while for $q=1$, 
\begin{equation*}		
a(n,k)=\begin{cases}
0, & 1\leq n<k; \\
1, & n=k; \\
n-k+1, & n>k. \end{cases}
\end{equation*}

We now formulate two results on necessary and sufficient conditions for the convergence of  series (\ref{BKseries}).

\begin{theorem}\label{theorem1}
	Let $-1\leq q<1$, $0<p<2$  and $r\geq p$.  Then  for any $\varepsilon>0$
	\begin{equation*}
	\sum_{n=1}^\infty n^{\frac{r}{p}-2}\mathbb{P} \Big\{\frac{|S_n|}{n^{1\!/\!p}}>\varepsilon  \Big\} <\infty \quad \text{iff} \quad \mathbb E|\theta|^{r}<\infty,
	\end{equation*}
	where $\mathbb E \theta=0 $ whenever $r\geq1$.
	
\end{theorem}

\begin{theorem}\label{theorem2}
	Let $q=1$,  $0<p<2/3$  and $r\geq p$.  Then for any  $ \varepsilon>0$
	\begin{equation*}
	\sum_{n=1}^\infty n^{\frac{r}{p}-2}\mathbb{P} \Big\{\frac{|S_n|}{n^{1\!/\!p}}>\varepsilon  \Big\} <\infty
	\quad \text{iff} \quad \mathbb E|\theta|^{\frac{r}{1-p}}<\infty,
	\end{equation*}
	where $\mathbb E \theta=0 $ whenever $r\geq1$.
	
\end{theorem}

The following example illustrates, that in general case for $p\geq 2/3$ the result of Theorem \ref{theorem2} is no longer true. 

\begin{example}
	Let $\theta$ be a standard normal random variable and 
	$\Phi_0 (x)=\frac{1}{\sqrt{2\pi}}\int_0^x e^{-x^2/2}\, dx$. Then the sum
	\begin{equation*}
	S_n=\sum_{k=1}^n (n-k+1)\theta_k, \ \ n\geq 1,
	\end{equation*}
	is normally distributed with zero mean and variance equal to $\dfrac{n}{6} (n+1)(2n+1)$.
	Therefore, a trivial calculation shows that 
	\begin{equation*}
	\mathbb{P}\Big\{|S_n|>n^{1/p}\varepsilon\Big\} =1-2\Phi_0\Big(\dfrac{n^{1/p}\varepsilon\sqrt{6}}{\sqrt{n(n+1)(2n+1)}} \Big)\underset{n\rightarrow\infty}{\centernot\longrightarrow}
	0,
	\end{equation*}
	if $p\geq 2/3$, which means that series {\rm{(\ref{BKseries})}} cannot converge whatever $r\geq p$.
\end{example}

\section{Proofs of main results}\label{proofs}

A crucial tool to prove both Theorem \ref{theorem1} and Theorem \ref{theorem2}, as in \cite{gut}, is a well-known Hoffmann-J$\o$rgensen inequality combined with moment inequalities for sums of r.v.'s. 
For ease of reading, we start this section recalling some probability inequalities used in the proofs below (see, for instance, \cite{LBai}).

Throughout the text, for a r.v. $\theta$ we denote by $\mu \theta$ (any of) its median and by $\theta^{'}$ an independent of $\theta$ and equidistributed with $\theta$ r. v., and set $\theta^{sym}=\theta-\theta^{'}$, i.e. $\theta^{sym}$ is a symmetrization of $\theta$.

\textbf{Weak symmetrization inequality.}\label{WS} 
Let $X$ be a r.v. For any $x$ and $a$  
\begin{gather*}
\dfrac12 \mathbb{P} \Big\{|X-\mu X| \geq x\Big\} \leq  \mathbb{P} \Big\{|X^{sym}| \geq x\Big\} \leq 2 \mathbb{P} \Big\{ |X-a|\geq x/2 \Big\}.
\end{gather*}

\textbf{Symmetrization moment inequality.}\label{SMI} 
Let $X$ be a r.v. For any $a$,
\begin{equation*}
\dfrac12 \mathbb{E} |X-\mu X|^m  \leq  \mathbb{E} |X^{sym}|^m \leq 2 c_r \mathbb{E} |X-a|^m, \ \ m>0,
\end{equation*}
where $c_r=1$ or $2^{r-1}$ in accordance with $r\leq 1$ or not.

\textbf{L\'{e}vy inequality.}\label{LV} 
Let $X_1$, $X_2$, ..., $X_n$ be independent symmetric r.v.'s and $S_n=\sum_{k=1}^n X_k$. Then for any  $x>0$,  
\begin{equation*}
\mathbb{P} \Big\{|S_n| >x\Big\} \geq \dfrac12 \mathbb{P} \Big\{\underset{1\leq j\leq n}{\max}|S_j|>x\Big\} \geq 
\dfrac12 \mathbb{P} \Big\{\underset{1\leq j\leq n}{\max}|X_j|>2x \Big\}.
\end{equation*}

\textbf{Hoffmann-J$\o$rgensen inequality for symmetric r.v.'s.}\label{HJ} 
Let $X_1$, $X_2$, ..., $X_n$ be independent symmetric r.v.'s,  $S_n=\sum_{k=1}^n X_k.$ Then for any $s,t >0$, 
\begin{equation}\label{HJ_s}
\mathbb{P} \Big\{|S_n| \geq 2t+s\Big\} \leq 
4\Big(\mathbb{P} \Big\{|S_n| \geq t\Big\}\Big)^2+
\mathbb{P} \Big\{\underset{1\leq j\leq n}{\max}|X_j| \geq s\Big\}.
\end{equation}

\textbf{Marcinkiewicz-Zygmund inequality.}\label{MZ} 
Let $r\geq 1$ and $X_1$, $X_2$, ..., $X_n$ be independent r.v.'s with $\mathbb{E} X_k =0$, $1\leq k\leq n$,  $S_n=\sum_{k=1}^n X_k$. Then there are positive constants $a_m \leq b_m$ such that  
\begin{equation*}
a_m \mathbb{E}\Big(\sum_{j=1}^n X_j^2 \Big)^{m/2} \leq \mathbb{E} |S_n|^m \leq
b_m \mathbb{E}\Big(\sum_{j=1}^n X_j^2 \Big)^{m/2}.
\end{equation*}

\textbf{$c_r$--inequality.}\label{CR} 
For r.v.'s $X_1$, $X_2$, ..., $X_n$ and  $S_n=\sum_{k=1}^n X_k$,  
\begin{equation*}
\mathbb{E} |S_n|^r \leq
c_r \sum_{j=1}^n \mathbb{E}|X_j|^r, 
\end{equation*}
where $c_r=1$ or $n^{r-1}$ according whether $0<r\leq 1$ or $r>1$.

Now proceed to the proof of Theorem \ref{theorem1}.

\renewcommand{\proofname}{Proof of Theorem \ref{theorem1}.}
\begin{proof}
	
	Let us start with the proof of sufficiency. First, note that due to uniform boundedness of the weights $a(n,k)$, for two partial cases, $r=p$  and $r=2p$, assertion of Theorem \ref{theorem1} is immediate from results by Gut, see \cite{gut}. Namely, convergence of ``Hsu-Robbins series'' (the case $r=2p$) and convergence of ``Spitzer series'' ($r=p$) follow from Theorem 7.1 and Theorem 7.4 in \cite{gut}, respectively.  
	Nevertheless, we carry out the proof for $r>p$ in all details.

	First restrict our considerations to the case of symmetrically distributed r.v. $\theta$.
	Let us fix any $\varepsilon>0$ and apply an iteration of Hoffmann-J$\o$rgensen inequality (\ref{HJ_s}) with $s=t=n^{1/p}\varepsilon$. Thus, for $j\geq1$ there exist some constants $C_j$ and $D_j$ such that 
	\begin{align}
	&\mathbb{P}\Big\{|S_n|>n^{1/p}\varepsilon\cdot 3^j\Big\} \leq \nonumber 
	\\
	&C_j \sum_{k=1}^n \mathbb{P}\Big\{ \Big|a(n,k) \theta_k \Big|>n^{1/p}\varepsilon\Big\}+D_j\Big(\mathbb{P}\Big\{|S_n|>n^{1/p}\varepsilon\Big\} \Big)^{2^j} =
	\nonumber 
	\\ 
	&C_j \sum_{k=1}^n \mathbb{P}\Big\{ \Big|\frac{1-q^{n-k+1}}{1-q}\theta_k \Big|>n^{1/p}\varepsilon\Big\}+D_j\Big(\mathbb{P}\Big\{|S_n|>n^{1/p}\varepsilon\Big\} \Big)^{2^j}.  \label{HJ1}
	\end{align}
	Note, that for $j=1$ we arrive at the classical version of Hoffmann-J$\o$rgensen inequality with $C_1=1$ and $D_1=4$.

	The first terms in (\ref{HJ1}) for $-1< q< 1$ can be estimated as follows
	\begin{align*}
	&\sum_{k=1}^n \mathbb{P} \Big\{\Big|\frac{1-q^{n-k+1}}{1-q} \theta_k\Big|>n^{1/p}\varepsilon\Big\}= 
	\\
	&\sum_{k=1}^n \mathbb{P} \Big\{(1-q^{n-k+1})|\theta_k|>n^{1/p}\varepsilon (1-q)\Big\}\leq
	\\
	&\sum_{k=1}^n \mathbb{P} \Big\{|\theta_k|>n^{1/p}\varepsilon_2\Big\}=n\mathbb{P} \Big\{|\theta|>n^{1/p}\varepsilon_2\Big\},
	\end{align*}
	where $\varepsilon_2= \varepsilon \cdot\big(1 \wedge (1-q)\big)$. For $q=-1$, however, 
	\begin{align}
	&\sum_{k=1}^n \mathbb{P} \Big\{\Big|\frac{1-q^{n-k+1}}{1-q} \theta_k\Big|>n^{1/p}\varepsilon\Big\}=
	\Big[\frac{n+1}{2} \Big] 
	\mathbb{P} \Big\{|\theta|>n^{1/p}\varepsilon \Big\}\leq \label{case q=-1} 
	\\
	&n\mathbb{P} \Big\{|\theta|>n^{1/p}\varepsilon\Big\}, \nonumber 
	\end{align}
	where $[\cdot]$ stands for the integer part.
	
	Without loss of generality set $\varepsilon_2=1$ and let $F_{\theta}(x)$ be the probability distribution function of $\theta$.
	Therefore, the first part of Baum-Katz series can be bounded as follows:
	\begin{align*}
	&\sum_{n=1}^\infty n^{r/p-2} \Big( \sum_{k=1}^n \mathbb{P} \Big\{\Big|\frac{1-q^{n-k+1}}{1-q} \theta_k\Big|>n^{1/p}\varepsilon\Big\}  \Big) \leq 
	\\
	&\sum_{n=1}^\infty n^{r/p-2} \cdot n \mathbb{P} \Big\{|\theta_k|>n^{1/p}\varepsilon_2\Big\}=
	2 \sum_{n=1}^\infty n^{r/p-1} \int_{n^{1/p}}^{\infty} dF_{\theta} (x)=
	\\
	&2 \int_{1}^{\infty}  \Big( \sum_{n=1}^{[x^p]} n^{r/p-1} \Big) dF_{\theta} (x) \sim 
	2 \int_{1}^{\infty} \Big( \int_1^{[x^p]} t^{r/p-1}\, dt \Big)dF_{\theta} (x)=
	\\
	&	2 \int_{1}^{\infty} \Big( \frac{p}{r} t^{r/p}\Big) \Big|_1^{[x^p]} dF_{\theta} (x)\sim
	\frac{2p}{r} \int_{1}^{\infty} x^r dF_{\theta} (x) <\infty, 
	\end{align*}
	since  $\mathbb E|\theta|^{r}<\infty$. By $ \sim $ we mean that both integrals are
	convergent or divergent simultaneously.

	Now switch to the second term in (\ref{HJ1}) and show that there exists some $j\geq 1$ that the series
	\begin{equation}\label{part of the series}
	\sum_{n=1}^\infty n^{r/p-2} \Big(\mathbb{P}\big\{|S_n|>n^{1/p}\varepsilon\big\}  \Big)^{2^j}
	\end{equation}
	converges. In order to do that let us find an upper bound
	for $\mathbb{P}\{|S_n|>n^{1/p}\varepsilon\}.$
	
	Firstly, by Markov inequality for $r>p$, we get
	\begin{equation*}
	\mathbb{P}\Big\{|S_n|>n^{1/p}\varepsilon\Big\}\leq \frac{\mathbb E|S_n|^r}{(n^{1/p}\varepsilon)^r}=\frac{\mathbb E|S_n|^r}{\varepsilon_1n^{r/p}}.
	\end{equation*}
	Next we deal with $\mathbb E|S_n|^r$ distinguishing between the following cases: 
	
	1) $0<r \leq 1$, 
	
	2) $r>1 $.
	
	\textbf{1)} Let $0<r \leq 1$. Applying $c_r$-inequality with $c_r=1$ to $\mathbb E|S_n|^r$, one obtains
	\begin{align*}
	&\mathbb E|S_n|^r \leq
	\sum_{k=1}^{n} \mathbb E\Big|\frac{1-q^{n-k+1}}{1-q} \theta_k \Big|^r=
	\\
	&\sum_{k=1}^{n} \Big(\frac{1-q^{n-k+1}}{1-q} \Big)^r \mathbb E |\theta_k|^r=
	\mathbb E|\theta|^r (1-q)^{-r} \sum_{k=1}^{n} \big(1-q^{n-k+1}\big)^r \leq
	\\
	&\mathbb E|\theta|^r (1-q)^{-r} \cdot 2^r n =
	C_1(r)  \mathbb E|\theta|^r n,
	\end{align*}
	where $C_1(r)=const=2^r(1-q)^{-r}.$

	\textbf{2)} Let $r>1$. In this case to $\mathbb E|S_n|^r$ we consequently apply  Marcinkiewicz-Zygmund inequality and the following well-known inequality: for positive $a_i$, $1\leq i \leq n$,  $n \in \mathbb N$ and $r>0$ it is true that 
	\begin{equation}\label{auxiliary inequality}
	(a_1^2+a_2^2+...+a_n^2)^{r/2}\leq n^{0\vee (r/2-1)} \sum_{i=1}^n a_i^{r}.
	\end{equation}
	Thus,
	\begin{align*}
	&\mathbb E|S_n|^r\leq 
	b_r \mathbb E\Big(\sum_{k=1}^{n} \Big(\frac{1-q^{n-k+1}}{1-q} \theta_k \Big)^2 \Big)^{r/2} \leq
	\\
	&b_r  n^{0\vee (r/2-1)}  \mathbb E \sum_{k=1}^{n} \Big(\frac{1-q^{n-k+1}}{1-q} \theta_k \Big)^{r} =
	\\
	&b_r  n^{0\vee (r/2-1)}  \mathbb E|\theta|^{r}(1-q)^{-r}  \sum_{k=1}^{n} \big(1-q^{n-k+1}\big)^r \leq  
	\\
	&b_r  n^{0\vee (r/2-1)}  \mathbb E|\theta|^{r} (1-q)^{-r} \cdot 2^r n= C_2(r)\mathbb E|\theta|^{r}  n^{1\vee (r/2)} ,
	\end{align*}
	where $C_2(r)=const=b_r2^r(1-q)^{-r}$.

	Let $C(r)$ denote a constant, which is equal to $C_1(r)$ or $C_2(r)$ depending on whether $r\leq1$ or $r>1$. Combining the above two cases, we arrive at the following bounds
	\begin{equation*}
	\mathbb E|S_n|^r \leq 
	C(r) \mathbb E|\theta|^r {n^{1 \vee (r/2)}}.
	\end{equation*}
	and 
	\begin{equation*}
	\mathbb{P}\Big\{|S_n|>n^{1/p}\varepsilon\Big\}\leq 
	\frac{\tilde{C}(r) \mathbb E|\theta|^r }{n^{r/p-(1 \vee (r/2))}},
	\end{equation*}
	where $\tilde{C}(r)=C(r)/\varepsilon_1$.
	
	Now, if $r\leq 2$, it is enough to set $j=1$ in (\ref{HJ1}) to obtain
	\begin{equation*}
	n^{r/p-2}  \Big(\mathbb{P}\big\{|S_n|>n^{1/p}\varepsilon\big\} \Big)^{2^j} \leq
	n^{r/p-2}  \frac{\big(C(r) \mathbb E|\theta|^{r} \big)^{2}} {n^{2r/p-2}}=\frac{\big(C(r) \mathbb E|\theta|^{r} \big)^{2}} { n^{r/p}},
	\end{equation*}
	and series (\ref{part of the series}) converges, since $\mathbb E|\theta|^{r}<\infty$.
	
	For $r>2$, 
	\begin{equation*}
	n^{r/p-2}  \Big(\mathbb{P}\big\{|S_n|>n^{1/p}\varepsilon\big\}  \Big)^{2^j} \leq
	n^{r/p-2}  \frac{\big(C(r) \mathbb E|\theta|^{r} \big)^{2^j}   } { n^{2^j(r/p-r/2)}}=\frac{\big(C(r) \mathbb E|\theta|^{r} \big)^{2}}{ n^{2^j\big(r/p-r/2\big)-r/p+2}} ,
	\end{equation*}
	whence series (\ref{part of the series})
	converges, provided that we choose $j$ so large that
	\begin{equation*}
	2^j \Big(\frac{r}{p}-\frac{r}{2}\Big)-\frac{r}{p}+2>1.
	\end{equation*}
	Thus, the proof of Theorem \ref{theorem1} for a symmetrically distributed r.v. $\theta$ is complete.

	Finally, show that sufficiency of Theorem \ref{theorem1} holds true for nonsymmetric $\theta$ as well. Indeed, according to symmetrization moment inequality, assumption $\mathbb E|\theta|^{r}<\infty$ implies that $\mathbb E|\theta^{sym}|^{r}<\infty$. The latter due to proved above yields   
	\begin{equation*}
	\sum_{n=1}^\infty n^{\frac{r}{p}-2}\mathbb{P} \Big\{|S_n^{sym}|>\varepsilon n^{1\!/\!p} \Big\}<\infty
	\end{equation*}
	with $S_n^{sym}$ being a symmetrization of $S_n$.  
	Notice that for any real $c$ the following set inclusion is true:
	\begin{equation*}
	\Big\{|S_n-S'_n|>c\Big\} \supset  \Big\{|S_n|>2c\Big\} \cap \Big\{|S'_n|\leq c\Big\},
	\end{equation*}
	whence
	\begin{equation*}
	\mathbb{P} \Big\{|S_n^{sym}|>c\Big\}\geq \mathbb{P} \Big\{|S_n|>2c\Big\}\cdot \mathbb{P} \Big\{|S'_n|\leq c\Big\},
	\end{equation*}
	where $S_n^{'}$ is an independent and equidistributed copy of $S_n$. The latter applied with $c=\varepsilon n^{1\!/\!p}$ yields
	\begin{equation*}
	\mathbb{P} \Big\{|S_n|>2\varepsilon n^{1\!/\!p}\Big\}  \leq \mathbb{P} \Big\{|S_n^{sym}|>\varepsilon n^{1\!/\!p}\Big\} \big/ \mathbb{P} \Big\{|S'_n|\leq \varepsilon n^{1\!/\!p}\Big\},
	\end{equation*}
	and one needs to have $\mathbb{P} \big\{|S'_n|\leq \varepsilon n^{1\!/\!p}\big\}$ bounded away from zero.
	But since $\mathbb E|\theta|^{r}<\infty$, with $\mathbb E \theta=0 $ whenever $r\geq1$, Weak Law of Large Numbers for rowwise independent random arrays (see Lemma 2.2, \cite{gut}) suggests that   
	$S'_n/n^{1/r}\underset{n\rightarrow \infty}{\longrightarrow} 0$ in probability. Finally, since $r\geq p$, then also $S'_n/n^{1/p}\underset{n\rightarrow \infty}{\longrightarrow} 0$ in probability. Therefore, $\mathbb{P} \big\{|S'_n|\leq \varepsilon n^{1\!/\!p}\big\}\underset{n\rightarrow \infty}{\longrightarrow} 1$.

	Now proceed to the proof of necessity. Although the methods used to prove this part are standard, let us first provide an idea of the proof in case of ``Hsu-Robbins series'', i.e. for $r=2p.$ 
	Initially assume that $\theta$ is symmetrically distributed and prove that   $\sum_{n=1}^\infty \mathbb{P} \big\{|S_n|>\varepsilon n^{1\!/\!p} \big\}< \infty$ implies  $ \mathbb E|\theta|^{2p}< \infty$. 
	
	According to representation \eqref{partial sums} and Theorem 2.3, \cite{gut}, the series   
	\begin{equation*}
	\sum_{n=1}^\infty \sum_{k=1}^n \mathbb{P} \Big\{|a(n,k)\theta_k|> \varepsilon n^{1\!/\!p} \Big \}
	\end{equation*}
	converges for any $\varepsilon>0$. Note that the latter conclusion is, in fact, corollary of Levy inequality combined with Borel-Cantelli Lemma. 
	Next, for $q=-1$ see equality (\ref{case q=-1}), while for $-1<q <1$,
	\begin{align*}
	&\sum_{k=1}^n \mathbb{P} \Big\{|a(n,k)\theta_k|> \varepsilon n^{1\!/\!p}\Big\}=
	\sum_{k=1}^n \mathbb{P} \Big\{\Big|\frac{1-q^{n-k+1}}{1-q}\theta_k\Big|> \varepsilon n^{\frac{1}{p}} \Big\}=
	\\
	&\sum_{k=1}^n \mathbb{P} \Big\{\Big|(1-q^{n-k+1})\theta_k\Big|> \varepsilon (1-q) n^{\frac{1}{p}} \Big\}\geq  n \mathbb{P} \Big\{ |\theta|>\varepsilon_1 n^{\frac{1}{p}}\Big\},
	\end{align*}
	where $\varepsilon_1=\varepsilon$, if $0\leq q<1$, and $\varepsilon_1=\frac{\varepsilon}{1+q}$, if $-1<q<0$.

	Let $F_{\theta}(x)$ be the probability distribution function of $\theta$. Without loss of generality, set $\varepsilon_1=1$, and consider the convergent series  
	\begin{align*}
	&\sum_{n=1}^\infty n \mathbb{P} \Big\{ |\theta|> n^{\frac{1}{p}}\Big\}=
	2 \sum_{n=1}^\infty n \int_{n^{\frac{1}{p}}}^\infty dF_{\theta}(x) =
	\\
	&2 \int_1^{\infty} \Big(\sum_{n=1}^{[x^p]} n \Big)dF_{\theta}(x)=
	2 \int_1^{\infty} \frac{(1+[x^p])[x^p]}{2} dF_{\theta}(x)> 
	\\
	&\int_1^{\infty} [x^p]^{2} \,dF_{\theta}(x)>
	\int_1^{\infty} \Big(\frac{x^p}{2}\Big)^2 dF_{\theta}(x)=
	\frac{1}{4}\int_1^{\infty} x^{2p}  \, dF_{\theta}(x).
	\end{align*}
	Therefore, the integral $\int_1^{\infty} x^{2p}   dF_{\theta}(x) $ converges, which yields $\mathbb E|\theta|^{2p}< \infty$.

	Now, let $\theta$ be an arbitrary random variable.
	According to weak symmetrization inequality $\sum_{n=1}^\infty \mathbb{P} \Big\{|S_n|>~\varepsilon n^{1\!/\!p} \Big\}<\infty$ implies $\sum_{n=1}^\infty \mathbb{P} \Big\{|S_n^{sym}|>\varepsilon n^{1\!/\!p}  \Big\}<\infty$. Therefore, in view of the above part of the proof one has $\mathbb E|\theta^{sym}|^{2p}< \infty$. Finally, again by symmetrization moment inequality,
	\begin{align*}
	&\mathbb E|\theta|^{2p} = \mathbb E|\theta- \mu\theta+\mu\theta|^{2p}\leq 
	C(p) \mathbb E \Big(|\theta- \mu\theta|^{2p}+|\mu\theta|^{2p}  \Big) \leq 
	\\
	&2C(p) \mathbb E|\theta^{sym}|^{2p}+C(p)|\mu\theta|^{2p},
	\end{align*}
	where $C(p)=1$ or $2^{2p-1}$ depending on whether $p\leq \frac{1}{2}$ or $\frac{1}{2} \leq p<2$. Thus, we proved that $\mathbb E|\theta|^{2p}<\infty$.

	Another important case is the ``Spitzer series'', that is when $r=p$. Here, to prove necessity, we follow the same idea with some additional reasoning, borrowed from \cite{gut}. Let us just provide some steps of it. 
	
	As before, assume that $\theta$ is symmetrically distributed and for any $\varepsilon>0$ the series $\sum_{n=1}^\infty \dfrac{1}{n} \mathbb{P} \Big\{|S_n|>\varepsilon n^{1\!/\!p}  \Big\}$ converges.
	Introduce a sequence of mutually independent events $(A_n, n\geq 1)$, which are independent with all the variables $\theta_k,$ $k\geq 1$, and $\mathbb{P} \{A_n\}=\frac{1}{n},$ $n\geq 1$. Then 
	\begin{equation*}
	\sum_{n=1}^\infty \frac{1}{n} \mathbb{P} \Big\{|S_n|>\varepsilon n^{1\!/\!p}  \Big\}=
	\sum_{n=1}^\infty  \mathbb{P} \Big\{A_n; |S_n|>\varepsilon n^{1\!/\!p} \Big\}.
	\end{equation*}
	Moreover, according to Levy inequality,
	\begin{equation*}
	\mathbb{P} \Big\{A_n; |S_n|>\varepsilon n^{1\!/\!p} \Big\}\geq \frac{1}{2} \mathbb{P} \Big\{A_n;\underset{1\leq k \leq n}\max|a(n,k)\theta_k|> 2\varepsilon n^{\frac{1}{p}} \Big\},
	\end{equation*}
	which yields  
	\begin{equation*}
	\sum_{n=1}^\infty \mathbb{P} \Big\{A_n;\underset{1\leq k \leq n}\max|a(n,k)\theta_k|> 2\varepsilon n^{\frac{1}{p}} \Big\}<\infty.
	\end{equation*}
	In view of Borel-Cantelli Lemma, 
	\begin{equation*}
	\mathbb{P} \Big\{A_n;\underset{1\leq k \leq n}\max|a(n,k)\theta_k|> 2\varepsilon n^{\frac{1}{p}} \ \ \textrm{infinitely often}  \Big\}=0.
	\end{equation*}
	The latter, in its turn, means that only finite number of rows, in which event $A_n$ takes place, have maximum terms greater than $2\varepsilon n^{\frac{1}{p}}$.
	Thus, if the rows were mutually independent, we could conclude according to Borel-Cantelli Lemma that 
	\begin{equation*}
	\sum_{n=1}^\infty \sum_{k=1}^n \mathbb{P} \Big\{A_n;|a(n,k)\theta_k|> 2\varepsilon n^{\frac{1}{p}}\Big\}<\infty,
	\end{equation*}
	which is what we need. But, since under the probability sign there are random variables, belonging to the same row, without loss of generality we may indeed assume the rows being independent, and therefore the series 
	\begin{equation*}
	\sum_{n=1}^\infty \frac{1}{n}\sum_{k=1}^n \mathbb{P} \Big\{|a(n,k)\theta_k|> \varepsilon n^{\frac{1}{p}} \Big\}
	\end{equation*}
	converges for any $\varepsilon>0$.
	To complete the proof of necessity in this case, one needs to literally follow the lines of it as in the previous case.

	By analogue, one can prove the necessity of Theorem \ref{theorem1} for any integer $r$, and then for any $r>p$. 
\end{proof}

\begin{remark}
	Due to constant coefficients in {\rm{\eqref{model}}} it became possible to prove the necessary part of Theorem~{\rm{\ref{theorem1}}} as well. Clearly, when coefficients $q_k$, $k\geq1$, depend on $k$, it is not as simple for both sufficiency and necessity.
	Nevertheless, for some specific cases of sequences with time-dependent coefficients, say, when $q_k$, $k\geq1$, are such that $|q_k|\leq q<1$, assumptions for the convergence of series {\rm{(\ref{BKseries})}} can be easily proved to be the same as in Theorem~{\rm{\ref{theorem1}}}.

\end{remark}

\renewcommand{\proofname}{Proof of Theorem \ref{theorem2}.}
\begin{proof}

	Sufficiency. First, note that in this case  
	\begin{equation*}
	\underset{1 \leq k\leq n} \sup |a(n,k)|\leq n, \ \ \text{for all} \ \ n\geq 1,
	\end{equation*}
	i.e. the weights $a(n,k)$ satisfy assumptions of Theorem 7.5, in \cite{gut} with $\lambda=1$ and $\alpha=1$. Therefore, sufficiency of Theorem \ref{theorem2} for $r=p$ immediately follows from Theorem 7.5, \cite{gut}. Note, that the case $r=2p$  is also immediate from Theorem 7.3, \cite{gut}. Nevertheless, we prove the sufficient part of Theorem \ref{theorem2} for $r>p$.

	In view of symmetrization-desymmetrization procedures, given in the proof of Theorem \ref{theorem1}, it suffices to restrict the proof to symmetrically distributed r.v. $\theta$. Moreover, since the proof of Theorem \ref{theorem2} follows exactly the same lines as the proof of Theorem \ref{theorem1}, we leave out some steps of it.

	Let us fix any $\varepsilon>0$ and apply an iteration of Hoffmann-J$\o$rgensen inequality (\ref{HJ_s}) with $s=t=n^{1/p}\varepsilon$. Thus, for $j\geq1$ there exist some constants $C_j$ and $D_j$ such that 
	\begin{align}
	&\mathbb{P}\Big\{|S_n|>n^{1/p}\varepsilon\cdot 3^j\Big\} \leq \nonumber \\
	&C_j \sum_{k=1}^n \mathbb{P}\Big\{ \Big|a(n,k) \theta_k \Big|>n^{1/p}\varepsilon\Big\}+D_j\Big(\mathbb{P}\Big\{|S_n|>n^{1/p}\varepsilon\Big\} \Big)^{2^j} =
	\nonumber \\ 
	&C_j \sum_{k=1}^n \mathbb{P}\Big\{ \Big|(n-k+1)\theta_k \Big|>n^{1/p}\varepsilon\Big\}+D_j\Big(\mathbb{P}\Big\{|S_n|>n^{1/p}\varepsilon\Big\} \Big)^{2^j}.  \label{HJ2}
	\end{align}

	The first term in (\ref{HJ2}) can be bounded as follows:
	\begin{equation*}
	\sum_{k=1}^n \mathbb{P} \Big\{\Big|(n-k+1) \theta_k\Big|>n^{1/p}\varepsilon\Big\}\leq 
	n\mathbb{P} \Big\{|\theta|>n^{1/p-1}\varepsilon\Big\}=2n \int_{n^{1/p-1}\varepsilon}^{\infty} dF_{\theta} (x),
	\end{equation*}
	where $F_{\theta}(x)$ is the probability distribution function of $\theta$.  Without loss of generality set $\varepsilon=1$.
	Then 
	\begin{align*}
	&\sum_{n=1}^\infty n^{r/p-2} \Big( \sum_{k=1}^n \mathbb{P} \Big\{\Big|(n-k+1) \theta_k\Big|>n^{1/p}\varepsilon\Big\}  \Big) \leq 
	\\
	&2\sum_{n=1}^\infty n^{r/p-2} \cdot n \int_{n^{1/p-1}}^{\infty} dF_{\theta}=
	2 \sum_{n=1}^\infty n^{r/p-1} \int_{n^{1/p}}^{\infty} dF_{\theta} (x)=
	\\
	&2 \int_{1}^{\infty}  \Big( \sum_{n=1}^{\big[x^{\frac{p}{1-p}}\big]} n^{r/p-1} \Big) dF_{\theta} (x) \sim 
	2 \int_{1}^{\infty} \Big( \int_1^{\big[x^{\frac{p}{1-p}}\big]} t^{r/p-1}\, dt \Big)\,dF_{\theta} (x)=
	\\
	&2 \int_{1}^{\infty} \Big( \frac{p}{r} t^{r/p}\Big) \Big|_1^{\big[x^{\frac{p}{1-p}}\big]} dF_{\theta} (x)\sim
	\frac{2p}{r} \int_{1}^{\infty} x^{\frac{r}{1-p}} \,dF_{\theta} (x) <\infty, 
	\end{align*}
	since  $\mathbb E|\theta|^{\frac{r}{1-p}}<\infty$.

	Now we deal with the second term in (\ref{HJ2}) and show that there exist some  $j \geq 1$ that 
	the series 
	\begin{equation}\label{part of the series 2}
	\sum_{n=1}^\infty n^{r/p-2} \Big(\mathbb{P}\big\{|S_n|>n^{1/p}\varepsilon\big\}  \Big)^{2^j}
	\end{equation}
	converges. To this end let us find some bounds 
	for $\mathbb{P}\{|S_n|>n^{1/p}\varepsilon\}$.
	
	Firstly, by Markov inequality,
	\begin{equation*}
	\mathbb{P}\Big\{|S_n|>n^{1/p}\varepsilon\Big\}\leq \frac{\mathbb E|S_n|^{r/(1-p)}}{(n^{1/p}\varepsilon)^{r/(1-p)}}=\frac{\mathbb E|S_n|^{r/(1-p)}}{\varepsilon_1 n^{\frac{r}{p}+\frac{r}{1-p}}}.
	\end{equation*}
	Next consider $\mathbb E|S_n|^{r/(1-p)}$, where $r>p$, distinguishing between such cases: 
	
	1) $ r/(1-p) \leq 1$, 
	
	2) $ r/(1-p)>1$.

	\textbf{1)} Let  $0<p< 1/2$ and $p<r\leq 1-p$. Applying $c_r$-inequality with $c_r=1$ to $\mathbb E|S_n|^{r/(1-p)}$, one obtains that 
	\begin{equation*}
	\mathbb E|S_n|^{\frac{r}{1-p}}\leq 
	\sum_{k=1}^{n} \mathbb E\Big|(n-k+1) \theta_k \Big|^{\frac{r}{1-p}}
	\leq\mathbb E|\theta|^{\frac{r}{1-p}}\cdot  n \cdot n^{\frac{r}{1-p}}=\mathbb E|\theta|^{\frac{r}{1-p}} \cdot n^{\frac{r}{1-p}+1}.
	\end{equation*}

	\textbf{2)} Let  $0<p<2/3$ and $r>p \vee  (1-p)$. In this case to $\mathbb E|S_n|^{r/(1-p)}$  we consequently apply Marcinkiewicz-Zygmund inequality and then inequality (\ref{auxiliary inequality}) with the power $r/(1-p)$ instead of $r$.
	Thus,
	\begin{align*}
	&\mathbb E|S_n|^{\frac{r}{1-p}}\leq
	b_r \mathbb E \Big(
	\sum_{k=1}^{n} \big(n-k+1\big)^2 |\theta_k|^2 \Big)^{\frac{r}{2(1-p)}}\leq 
	\\
	&b_r \cdot  n^{0\vee \big(\frac{r}{2(1-p)}-1\big)}  
	\sum_{k=1}^{n} \big(n-k+1\big)^{\frac{r}{(1-p)}} \mathbb E|\theta_k|^{\frac{r}{(1-p)}} =
	\\
	&b_r \cdot n^{0\vee \big(\frac{r}{2(1-p)}-1\big)}   \mathbb E|\theta|^{\frac{r}{1-p}} \cdot  \sum_{k=1}^{n} \big(n-k+1\big)^{\frac{r}{(1-p)}} \leq
	\\
	&b_r \cdot n^{0\vee \big(\frac{r}{2(1-p)}-1\big)}   \mathbb E|\theta|^{\frac{r}{1-p}} \cdot n \cdot n^{\frac{r}{1-p}}=
	b_r \mathbb E|\theta|^{\frac{r}{1-p}}\cdot n^{\frac{r}{1-p}+\big(1 \vee \frac{r}{2(1-p)}\big)},
	\end{align*}
	where $b_r$ is some constant depending on r.

	Combining together cases 1) and 2), we get the following bounds:
	\begin{equation*}
	\mathbb E|S_n|^r \leq 
	B(r) \mathbb E|\theta|^{\frac{r}{1-p}}\cdot n^{\frac{r}{1-p}+\big(1 \vee \frac{r}{2(1-p)}\big)}.
	\end{equation*}
	and 
	\begin{equation*}
	\mathbb{P}\Big\{|S_n|>n^{1/p}\varepsilon\Big\}\leq 
	\frac{\tilde{B}(r) \mathbb E|\theta|^{\frac{r}{1-p}}\cdot n^{\frac{r}{1-p}+\big(1 \vee \frac{r}{2(1-p)}\big)}}
	{n^{\frac{r}{p}+\frac{r}{1-p}}}=
	\frac{\tilde{B}(r) \mathbb E|\theta|^{\frac{r}{1-p}}}
	{n^{\frac{r}{p}-\big(1 \vee \frac{r}{2(1-p)}\big)}},
	\end{equation*}
	where $\tilde{B}(r)=B(r)/\varepsilon_1$, and $B(r)=1$ or $b_r$ according whether 
	$ r/(1-p) \leq 1$ or $ r/(1-p)>1$. Hence
	\begin{equation*}
	n^{r/p-2}\cdot \Big(\mathbb{P}\big\{|S_n|>n^{1/p}\varepsilon\big\}  \Big)^{2^j}\leq 
	\frac{\tilde{B}(r)^{2^j} \Big(\mathbb E|\theta|^{\frac{r}{1-p}} \Big)^{2^j} n^{r/p-2}}{n^{2^j\Big(\frac{r}{p}-\big(1 \vee \frac{r}{2(1-p)}\big) \Big)}}.
	\end{equation*}
	
	Finally, if $\frac{r}{(1-p)} \leq 2$ it suffices to set $j=1$ in (\ref{HJ2}), whence
	\begin{equation*}
	n^{r/p-2}  \Big(\mathbb{P}\big\{|S_n|>n^{1/p}\varepsilon\big\}  \Big)^{2^j} \leq
	\frac{\Big(\tilde{B}(r) \mathbb E|\theta|^{\frac{r}{1-p}} \Big)^{2}}{n^{r/p}},
	\end{equation*}
	and the series (\ref{part of the series 2})
	converges, since $\mathbb E|\theta|^{\frac{r}{1-p}}<\infty$ and $r>p$.

	When $\frac{r}{(1-p)} > 2$ the series (\ref{part of the series 2})  is convergent, if $\mathbb E|\theta|^{\frac{r}{1-p}}<\infty$ and
	\begin{equation*}
	2^j \Big(\frac{r}{p}-\frac{r}{2(1-p)}\Big)-\frac{r}{p}+2>1.
	\end{equation*}
	But, according to assumptions imposed on $p$ and $r$, it is always possible to pick $j$ so large that the latter inequality holds true.

	We avoid repeating the necessary part of the proof for it is fully based on the same ideas as in Theorem \ref{theorem1}. 
\end{proof}

\end{document}